\documentclass[a4paper]{amsart}
\usepackage{amssymb} 
\usepackage[T1]{fontenc}
\usepackage[latin1]{inputenc}
\usepackage{amsfonts}
\usepackage{amsxtra}
\usepackage{ae}
\usepackage[all]{xy}
\usepackage{enumerate}

\include{diagram}

\newcommand*{\ket}{\rangle}
\newcommand*{\bra}{\langle}

\newcommand*{\D}{\mathcal{D}}

\newcommand*{\cotimes}{\hat{\otimes}}
\newcommand*{\ayd}{\mathsf{AYD}}
\newcommand*{\AYD}{\mathsf{A}}
\newcommand*{\hit}{\rightharpoonup}
\newcommand*{\hitby}{\leftharpoonup}

\DeclareMathOperator{\LSMod}{-\mathsf{Mod}}
\DeclareMathOperator{\SComodR}{\mathsf{Comod}-}
\DeclareMathOperator{\pro}{pro}

\DeclareMathOperator{\Hom}{Hom}

\DeclareMathOperator{\id}{id}

\newenvironment{bnum}
{\begin{list}{}
    {\setlength{\labelwidth}{15pt}
     \setlength{\leftmargin}{\labelwidth}
    }
}
{\end{list}}

\numberwithin{equation}{section}
\theoremstyle{change}
\newtheorem{theorem}{Theorem}[section]
\newtheorem{prop}[theorem]{Proposition}
\newtheorem{lemma}[theorem]{Lemma}

\newtheorem{definition}[theorem]{Definition}

\begin{document}

\title[Equivariant cyclic homology for quantum groups]{Equivariant cyclic homology for 
quantum groups}
\author{Christian Voigt}
\address{Institut for Mathematical Sciences\\
         University of Copenhagen\\
         Universitetsparken 5 \\
         2100 Copenhagen\\
         Denmark
}
\email{cvoigt@math.ku.dk}


\maketitle

\begin{abstract}
We define equivariant periodic cyclic homology for bornological quantum groups. Generalizing corresponding 
results from the group case, we show that the theory is homotopy invariant, stable and satisfies excision in 
both variables. Along the way we prove Radfords formula for the antipode of a bornological quantum group. 
Moreover we discuss anti-Yetter-Drinfeld modules and 
establish an analogue of the Takesaki-Takai duality theorem in the setting of bornological quantum groups. 
\end{abstract}

\section{Introduction}

Equivariant cyclic homology can be viewed as a noncommutative generalization of equivariant de Rham cohomology. 
For actions of finite groups or compact Lie groups, different aspects of the theory have been studied by various 
authors \cite{BG}, \cite{Brylinski1}, \cite{Brylinski2}, \cite{KKL1}, \cite{KKL2}. 
In order to treat noncompact groups as well, a general framework for equivariant cyclic homology following the 
Cuntz-Quillen formalism \cite{CQ1}, \cite{CQ2}, \cite{CQ3} has been introduced in \cite{Voigtepch}. For instance, in the setting of 
discrete groups or totally disconnected groups this yields a new approach to classical 
constructions in algebraic topology \cite{Voigtbs}. However, in contrast 
to the previous work mentioned above, a crucial feature of the construction in \cite{Voigtepch} is the fact that the 
basic ingredient in the theory is not a complex in the usual sense of homological algebra. In particular, 
the theory does not fit into the traditional scheme of defining cyclic homology using cyclic modules or 
mixed complexes. \\
In this note we define equivariant periodic cyclic homology for quantum groups. 
This generalizes the constructions in the group case developped in \cite{Voigtepch}. Again we work 
in the setting of bornological vector spaces. Correspondingly, the appropriate notion 
of a quantum group in this context is the concept of a bornological quantum group introduced in \cite{Voigtbqg}. 
This class of quantum groups includes all locally compact groups and their duals as well 
as all algebraic quantum groups in the sense of van Daele \cite{vD2}. As in the theory of van Daele, an important ingredient 
in the definition of a bornological quantum group is the Haar measure. It is crucial for the duality theory and 
also explicitly used at several points in the construction of the homology theory presented in this paper. However, with some modifications 
our definition of equivariant cyclic homology could also be adapted to a completely algebraic setting using Hopf algebras 
with invertible antipodes instead. \\
From a conceptual point of view, equivariant cyclic homology should be viewed as a homological analogon to equivariant 
$ KK $-theory \cite{Kasparov1}, \cite{Kasparov2}. The latter has been extended by Baaj and Skandalis to  
coactions of Hopf-$ C^* $-algebras \cite{BaSka1}. However, in our situation it is more convenient to work with 
actions instead of coactions. \\
An important ingredient in equivariant cyclic homology is the concept of a covariant module \cite{Voigtepch}. 
In the present paper we will follow the more common terminology introduced in \cite{HKRS} and call these 
objects anti-Yetter-Drinfeld modules instead. In order to construct the natural symmetry operator on these modules in the 
general quantum group setting we prove 
a formula relating the fourth power of the antipode with the modular functions of a bornological quantum group and its dual. 
In the context of finite dimensional Hopf algebras this formula is a classical result due to Radford \cite{Radford}. \\
Although anti-Yetter-Drinfeld modules occur naturally in the constructions one should point out that our theory does not fit into 
the framework of Hopf-cyclic cohomology \cite{HKRScyclic}. Still, there are relations to previous 
constructions for Hopf algebras by Akbarpour and Khalkhali \cite{BK1}, \cite{BK2} as well as Neshveyev and Tuset \cite{NT}. 
Remark in particular that cosemisimple Hopf algebras or finite dimensional Hopf algebras can be viewed as bornological quantum groups as well. 
However, basic examples show that 
the homology groups defined in \cite{BK1}, \cite{BK2}, \cite{NT} only reflect a small part of the information contained 
in the theory described below. \\
Let us now describe how the paper is organized. In section \ref{secqg} we 
recall the definition of a bornological quantum group. We explain some basic features 
of the theory including the definition of the dual quantum group and the Pontrjagin duality theorem. 
This is continued in section \ref{secaction} where we 
discuss essential modules and comodules over bornological quantum groups as well as actions on 
algebras and their associated crossed products. We prove an analogue of the Takesaki-Takai duality theorem 
in this setting. Section \ref{secradford} contains the discussion of Radford's formula relating the 
antipode with the modular functions of a quantum group and its dual. 
In section \ref{seccovmod} we study anti-Yetter-Drinfeld modules over bornological quantum groups and 
introduce the notion of a paracomplex. 
After these preparations we define equivariant periodic cyclic homology in section \ref{secech}. Finally, in 
section \ref{sechomstabex} we show that our theory is homotopy invariant, stable 
and satisfies excision in both variables. \\
Throughout the paper we work over the complex numbers. For simplicity we have avoided 
the use of pro-categories in connection with the Cuntz-Quillen 
formalism to a large extent. 

\section{Bornological quantum groups} \label{secqg} 

The notion of a bornological quantum group was introduced in \cite{Voigtbqg}. We will work with this concept  
in our approach to equivariant cyclic homology. For information on bornological vector spaces and more details 
we refer to \cite{H-L2}, \cite{Meyerthesis}, \cite{Voigtbqg}. All bornological vector spaces are 
assumed to be convex and complete. \\
A bornological algebra $ H $ is called essential if the multiplication map induces an isomorphism $ H \cotimes_H H \cong H $. 
The multiplier algebra $ M(H) $ of a bornological algebra $ H $ consists of all two-sided multipliers of $ H $, the latter being 
defined by the usual algebraic conditions. There exists a canonical bounded 
homomorphism $ \iota: H \rightarrow M(H) $. A bounded linear functional $ \phi: H \rightarrow \mathbb{C} $ on a bornological algebra 
is called faithful if $ \phi(xy) = 0 $ for all $ y \in H $ implies $ x = 0 $ and $ \phi(xy) = 0 $ for all $ x $ implies $ y = 0 $. 
If there exists such a functional the map $ \iota: H \rightarrow M(H) $ is injective. 
In this case one may view $ H $ as a subset of the multiplier algebra $ M(H) $. \\
In the sequel $ H $ will be an essential bornological algebra with a faithful bounded linear functional. 
For technical reasons we assume moreover that the underlying bornological vector space of $ H $ satisfies the approximation 
property. \\
A module $ M $ over $ H $ is called essential if the module action induces an isomorphism $ H \cotimes_H M \cong M $. 
Moreover an algebra homomorphism $ f: H \rightarrow M(K) $ is essential if $ f $ turns $ K $ into an essential left and right 
module over $ H $. Assume that $ \Delta: H \rightarrow M(H \cotimes H) $ is an essential homomorphism. 
The map $ \Delta $ is called a comultiplication if it is coassociative, that is, if $ (\Delta \cotimes \id) \Delta = 
(\id \cotimes \Delta)\Delta $ holds. Moreover the Galois maps $ \gamma_l, \gamma_r, \rho_l, \rho_r : H \cotimes H \rightarrow M(H \cotimes H) $ 
for $ \Delta $ are defined by 
\begin{align*}
\gamma_l(x\otimes y) &= \Delta(x)(y \otimes 1), &\gamma_r(x \otimes y) = \Delta(x) (1 \otimes y) \\ 
\rho_l(f \otimes g) &= (x \otimes 1) \Delta(y), &\rho_r(x \otimes y) = (1 \otimes x) \Delta(y). 
\end{align*}
Let $ \Delta: H \rightarrow M(H \cotimes H) $ be a comultiplication such that all Galois maps associated 
to $ \Delta $ define bounded linear maps from $ H \cotimes H $ into itself. If $ \omega $ is a bounded linear 
functional on $ H $ we define for every $ x \in H $ a multiplier $ (\id \cotimes \omega)\Delta(x) \in M(H) $ by 
\begin{eqnarray*}
(\id \cotimes \omega)\Delta(x) \cdot y = (\id \cotimes \omega)\gamma_l(x \otimes y) \\
y \cdot (\id \cotimes \omega)\Delta(x) = (\id \cotimes \omega)\rho_l(y \otimes x). 
\end{eqnarray*}
In a similar way one defines $ (\omega \cotimes \id)\Delta(x) \in M(H) $.
A bounded linear functional $ \phi: H \rightarrow \mathbb{C} $ is called left invariant if 
\begin{equation*}
(\id \cotimes \phi)\Delta(x) = \phi(x) 1
\end{equation*}
for all $ x \in H $. Analogously one defines right invariant functionals. \\
Let us now recall the definition of a bornological quantum group.
\begin{definition}\label{bqgdef}
A bornological quantum group is an essential bornological algebra $ H $ satisfying the approximation property 
with a comultiplication $ \Delta: H \rightarrow M(H \cotimes H) $ such that all 
Galois maps associated to $ \Delta $ are isomorphisms together with a faithful left invariant 
functional $ \phi: H \rightarrow \mathbb{C} $. 
\end{definition}
The definition of a bornological quantum group is equivalent to the definition of an algebraic quantum group in the 
sense of van Daele \cite{vD2} provided the underlying bornological vector space carries the fine bornology. The functional 
$ \phi $ is unique up to a scalar and referred to as the left Haar functional of $ H $. 
\begin{theorem} \label{bqchar}
Let $ H $ be a bornological quantum group. Then there exists an essential algebra homomorphism 
$ \epsilon: H \rightarrow \mathbb{C} $ and a linear isomorphism $ S: H \rightarrow H $ which is both an 
algebra antihomomorphism and a coalgebra antihomomorphism such that 
\begin{equation*}
(\epsilon \cotimes \id)\Delta = \id = (\id \cotimes \epsilon)\Delta
\end{equation*}
and
\begin{equation*}
\mu(S \cotimes \id) \gamma_r = \epsilon \cotimes \id, \qquad \mu(\id \cotimes S) \rho_l = \id \cotimes \epsilon.
\end{equation*}
Moreover the maps $ \epsilon $ and $ S $ are uniquely determined. 
\end{theorem}
Using the antipode $ S $ one obtains that every bornological quantum group is equipped with a faithful right 
invariant functional $ \psi $ as well. Again, such a functional is uniquely determined up to a scalar. There are
injective bounded linear maps $ \mathcal{F}_l, \mathcal{F}_r, \mathcal{G}_l, \mathcal{G}_r: 
H \rightarrow H' = \Hom(H, \mathbb{C}) $ defined by the formulas
\begin{align*}
\mathcal{F}_l(x)(h) &= \phi(hx), \qquad \mathcal{F}_r(x)(h) = \phi(xh) \\
\mathcal{G}_l(x)(h) &= \psi(hx), \qquad \mathcal{G}_r(x)(h) = \psi(xh).
\end{align*} 
The images of these maps coincide and determine a vector space $ \hat{H} $. Moreover, there exists a unique bornology 
on $ \hat{H} $ such that these maps are bornological isomorphisms. 
The bornological vector space $ \hat{H} $ is equipped with a multiplication which is induced from the comultiplication of 
$ H $. In this way $ \hat{H} $ becomes an essential bornological algebra and the multiplication of 
$ H $ determines a comultiplication on $ \hat{H} $. 
\begin{theorem} 
Let $ H $ be a bornological quantum group. Then $ \hat{H} $ with the structure maps described above is again a 
bornological quantum group. The dual quantum group of $ \hat{H} $ is canonically isomorphic to $ H $.
\end{theorem}
Explicitly, the duality isomorphism $ P: H \rightarrow \hat{\hat{H}} $ is given by $ P = \hat{\mathcal{G}}_l \mathcal{F}_l S $ 
or equivalently $ P = \hat{\mathcal{F}}_r \mathcal{G}_r S $. Here we write $ \hat{\mathcal{G}}_l $ and $ \hat{\mathcal{F}}_r $ 
for the maps defined above associated to the dual Haar functionals on $ \hat{H} $. The second statement of the previous theorem should be viewed as an 
analogue of the Pontrjagin duality theorem. \\
In \cite{Voigtbqg} all calculations were written down explicitly in terms of the Galois maps and their inverses. However, in this way many 
arguments tend to become lengthy and not particularly transparent. To avoid this we shall use the Sweedler notation in the 
sequel. That is, we write
$$
\Delta(x) = x_{(1)} \otimes x_{(2)}
$$
for the coproduct of an element $ x $, and accordingly for higher coproducts. Of course this has to be handled with care since 
expressions like the previous one only have a formal meaning. Firstly, the element $ \Delta(x) $ is a multiplier and not contained in an actual 
tensor product. Secondly, we work with completed tensor products which means that even a generic element in $ H \cotimes H $ cannot be 
written as a finite sum of elementary tensors as in the algebraic case. 

\section{Actions, coactions and crossed products} \label{secaction}

In this section we review the definition of essential comodules over a bornological quantum group and their 
relation to essential modules over the dual. 
Moreover we consider actions on algebras and their associated crossed products and prove an analogue of the Takesaki-Takai 
duality theorem. \\
Let $ H $ be a bornological quantum group. Recall from section \ref{secqg} that a module $ V $ over $ H $ is called essential if 
the module action induces an isomorphism $ H \cotimes_H V \cong V $. A bounded linear map $ f: V \rightarrow W $ between 
essential $ H $-modules is called $ H $-linear or $ H $-equivariant if it commutes with the module actions. We denote the 
category of essential $ H $-modules and equivariant linear maps by $ H \LSMod $. Using the comultiplication of $ H $ one obtains 
a natural $ H $-module structure on the tensor product of two $ H $-modules and $ H \LSMod $ becomes a monoidal category in this way. \\
We will frequently use the regular actions associated to a bornological quantum group $ H $. 
For $ t \in H $ and $ f \in \hat{H} $ one defines 
$$
t \hit f = f_{(1)} \, f_{(2)}(t), \qquad f \hitby t = f_{(1)}(t) f_{(2)}
$$
and this yields essential left and right $ H $-module structures on $ \hat{H} $, respectively. \\
Dually to the concept of an essential module one has the notion of an essential comodule. Let $ H $ be a bornological quantum 
group and let $ V $ be a bornological vector space. 
A coaction of $ H $ on $ V $ is a right $ H $-linear bornological isomorphism $ \eta: V \cotimes H \rightarrow V \cotimes H $ 
such that the relation
$$
(\id \otimes \gamma_r) \eta_{12} (\id \otimes \gamma_r^{-1}) = \eta_{12} \eta_{13}
$$
holds. 
\begin{definition}
Let $ H $ be a bornological quantum group. An essential $ H $-comodule is a bornological vector space 
$ V $ together with a coaction $ \eta: V \cotimes H \rightarrow V \cotimes H $. 
\end{definition}
A bounded linear map $ f: V \rightarrow W $ between essential comodules is called $ H $-colinear if 
it is compatible with the coactions in the obvious sense. We write $ \SComodR H $ for the category of 
essential comodules over $ H $ with $ H $-colinear maps as morphisms. The category $ \SComodR H $ is a monoidal 
category as well. \\
If the quantum group $ H $ is unital, a coaction is the same thing as 
a bounded linear map $ \eta: V \rightarrow  V \cotimes H $ such that 
$ (\eta \cotimes \id) \eta = (\id \cotimes \Delta)\eta $ and $ (\id \cotimes \epsilon) \eta = \id $. \\
Modules and comodules over bornological quantum groups are related in the same way as modules 
and comodules over finite dimensional Hopf algebras. 
\begin{theorem} \label{modcomod}
Let $ H $ be a bornological quantum group. 
Every essential left $ H $-module is an essential right $ \hat{H} $-comodule in a natural way and vice versa.  
This yields inverse isomorphisms between the category of essential $ H $-modules and the category 
of essential $ \hat{H} $-comodules. These isomorphisms are compatible with tensor products. 
\end{theorem}
Since it is more convenient to work with essential modules instead of comodules
we will usually prefer to consider modules in the sequel. \\
An essential $ H $-module is called projective if it has the lifting property with respect to surjections of essential $ H $-modules 
with bounded linear splitting. It is shown in \cite{Voigtbqg} that a bornological quantum group $ H $ is projective as a left 
module over itself. This can be generalized as follows. 
\begin{lemma} \label{omegaprojective}
Let $ H $ be a bornological quantum group and let $ V $ be any essential $ H $-module. Then the essential $ H $-modules 
$ H \cotimes V $ and $ V \cotimes H $ are projective.
\end{lemma}
\proof Let $ V_\tau $ be the space $ V $ equipped with the trivial $ H $-action induced by the counit. 
We have a natural $ H $-linear isomorphism $ \alpha_l: H \cotimes V \rightarrow H \cotimes V_\tau $ 
given by $ \alpha_l(x \otimes v) = x_{(1)} \otimes S(x_{(2)}) \cdot v $. Similarly, the map 
$ \alpha_r: V \cotimes H \rightarrow V_\tau \cotimes H $ given by $ \alpha_r(v \otimes x) = S^{-1}(x_{(1)}) \cdot v \otimes x_{(2)} $ 
is an $ H $-linear isomorphism. Since $ H $ is projective this yields the claim. \qed \\
Using category language an $ H $-algebra is by definition an algebra in the category $ H \LSMod $. We formulate this more explicitly in the 
following definition.
\begin{definition}
Let $ H $ be a bornological quantum group. An $ H $-algebra is a bornological algebra $ A $ which is at the same 
time an essential $ H $-module such that the multiplication map $ A \cotimes A \rightarrow A $ is $ H $-linear.
\end{definition}
If $ A $ is an $ H $-algebra we will also speak of an action of $ H $ on $ A $. Remark that we do not assume that an algebra has 
an identity element. The unitarization $ A^+ $ of an $ H $-algebra $ A $ becomes an $ H $-algebra by considering 
the trivial action on the extra copy $ \mathbb{C} $. \\
According to theorem \ref{modcomod} we can equivalenty describe an $ H $-algebra as a bornological algebra 
$ A $ which is at the same time an essential $ \hat{H} $-comodule such that the multiplication is $ \hat{H} $-colinear. 
Under additional assumptions there is another possibility to describe this structure which resembles the definition 
of a coaction in the setting of $ C^* $-algebras. 
\begin{definition}
Let $ H $ be a bornological quantum group. An algebra coaction of $ H $ on an essential bornological algebra $ A $ 
is an essential algebra homomorphism $ \alpha: A \rightarrow M(A \cotimes H) $ such that the maps 
$ \alpha_l $ and $ \alpha_r $ from $ A \cotimes H $ to $ M(A \cotimes H) $ given by 
$$
\alpha_l(a \otimes x) = (1 \otimes x)\alpha(a), \qquad \alpha_r(a \otimes x) = \alpha(a)(1 \otimes x)
$$
induce bornological automorphisms of $ A \cotimes H $ and 
$$
(\alpha \cotimes \id) \alpha = (\id \cotimes \Delta)\alpha.  
$$
\end{definition}
Let us call an essential bornological algebra $ A $ regular if it is equipped with a faithful bounded linear functional 
and satisfies the approximation property. An algebra coaction $ \alpha: A \rightarrow M(A \cotimes H) $ on a regular 
bornological algebra $ A $ is injective and one has $ (\id \cotimes \epsilon) \alpha = \id $. 
\begin{prop}
Let $ H $ be a bornological quantum group and let $ A $ be a regular bornological algebra. 
Then every algebra coaction of $ \hat{H} $ on $ A $ corresponds to a unique
$ H $-algebra structure on $ A $ and vice versa. 
\end{prop}
\proof Assume that $ \alpha $ is an algebra coaction of $ \hat{H} $ on $ A $ and define $ \eta = \alpha_r $. By definition 
$ \eta $ is a right $ \hat{H} $-linear automorphism of $ A \cotimes \hat{H} $. We compute
\begin{align*}
(\id \cotimes \gamma_r)&\eta^{12}(\id \cotimes \gamma_r^{-1})(a \otimes f \otimes g) = 
(\id \cotimes \gamma_r)(\alpha(a)(1 \otimes \gamma_r^{-1}(f \otimes g))) \\
&= (\id \cotimes \Delta)(\alpha(a))(1 \otimes \gamma_r \gamma_r^{-1}(f \otimes g)) \\
&= (\alpha \cotimes \id)(\alpha(a))(1 \otimes f \otimes g) \\
&= \eta^{12} \eta^{13}(a \otimes f \otimes g)
\end{align*} 
which shows that $ \eta $ defines a right $ \hat{H} $-comodule structure on $ A $. Moreover we have 
\begin{align*}
(\mu \cotimes \id)\eta^{13}&\eta^{23}(a \otimes b \otimes f) = (\mu \cotimes \id)\eta^{13}(a \otimes \alpha(b)(1 \otimes f)) \\
&= \alpha(a) \alpha(b)(1 \otimes f) = \alpha(ab)(1 \otimes f) = \eta(\mu \cotimes \id)(a \otimes b \otimes f)
\end{align*} 
and it follows that $ A $ becomes an $ H $-algebra using this coaction. \\
Conversely, assume that $ A $ is an $ H $-algebra implemented by the coaction $ \eta: A \cotimes \hat{H} \rightarrow A \cotimes \hat{H} $. 
Define bornological automorphisms $ \eta_l $ and $ \eta_r $ of $ A \cotimes \hat{H} $ by 
$$
\eta_l = (\id \cotimes S^{-1})\eta^{-1}(\id \cotimes S), \qquad \eta_r = \eta. 
$$
The map $ \eta_l $ is left $ \hat{H} $-linear for the action of $ \hat{H} $ on the second tensor factor and 
$ \eta_r $ is right $ \hat{H} $-linear. 
Since $ \eta $ is a compatible with the multiplication we have 
$$ 
\eta_r(\mu \cotimes \id) = (\mu \cotimes \id) \eta_r^{13} \eta_r^{23} 
$$
and  
$$
\eta_l(\mu \cotimes \id) = (\mu \cotimes \id) \eta_l^{23} \eta_l^{13}.
$$
In addition one has the equation
$$
(\id \cotimes \mu)\eta_l^{12} = (\id \cotimes \mu) \eta_r^{13} 
$$
relating $ \eta_l $ and $ \eta_r $. These properties of the maps $ \eta_l $ and $ \eta_r $ imply that
$$
\alpha(a)(b \otimes f) = \eta_r(a \otimes f)(b \otimes 1), \qquad 
(b \otimes f) \alpha(a) = (b \otimes 1) \eta_l(a \otimes f)
$$
defines an algebra homomorphism $ \alpha $ from $ A $ to $ M(A \cotimes \hat{H}) $. 
As in the proof of proposition 7.3 in \cite{Voigtbqg} one shows that $ \alpha $ is essential. 
Observe that we may identify the natural map $ A \cotimes_A (A \cotimes \hat{H}) \rightarrow A \cotimes \hat{H} $ 
induced by $ \alpha $ with $ \eta_r^{13}: A \cotimes_A (A \cotimes \hat{H} \cotimes_{\hat{H}} \hat{H}) \rightarrow 
(A \cotimes_A A) \cotimes (\hat{H} \cotimes_{\hat{H}} \hat{H}) $ 
since $ A $ is essential. \\
The maps $ \alpha_l $ and $ \alpha_r $ associated to the homomorphism $ \alpha $ can be identified 
with $ \eta_l $ and $ \eta_r $, respectively. Finally, the coaction identity 
$ (\id \cotimes \gamma_r)\eta^{12}(\id \cotimes \gamma_r^{-1}) = \eta^{12} \eta^{13} $ implies 
$ (\alpha \cotimes \id)\alpha = (\id \cotimes \Delta)\alpha $. Hence $ \alpha $ defines an algebra coaction of 
$ \hat{H} $ on $ A $. \\
It follows immediately from the constructions that the two procedures described above are inverse to each other. \qed \\
To every $ H $-algebra $ A $ one may form the associated crossed product $ A \rtimes H $. The underlying bornological 
vector space of $ A \rtimes H $ is $ A \cotimes H $ and the multiplication is defined by the chain of maps
$$
 \xymatrix{
     A \cotimes H \cotimes A \cotimes H \; \ar@{->}[r]^{\gamma_r^{24}} & 
     A \cotimes H \cotimes A \cotimes H \; \ar@{->}[r]^{\quad \id \cotimes \lambda \cotimes \id} & 
     \; A \cotimes A \cotimes H \ar@{->}[r]^{\;\;\;\mu \cotimes \id} & A \cotimes H
     }
$$
where $ \lambda $ denotes the action of $ H $ on $ A $. Explicitly, the multiplication in $ A \rtimes H $ is given by the formula
\begin{equation*}
(a \rtimes x)(b \rtimes y) = a x_{(1)} \cdot b \otimes x_{(2)}y
\end{equation*}
for $ a, b \in A $ and $ x, y \in H $. On the crossed product $ A \rtimes H $ one has the dual action of $ \hat{H} $ defined by 
$$
f \cdot (a \rtimes x) = a \rtimes (f \hit x)
$$
for all $ f \in \hat{H} $. In this way $ A \rtimes H $ becomes an $ \hat{H} $-algebra. Consequently 
one may form the double crossed product $ A \rtimes H \rtimes \hat{H} $. In the remaing part of 
this section we discuss the Takesaki-Takai duality isomorphism which clarifies the structure 
of this algebra. \\
First we describe a general construction which will also be needed later in connection with stability 
of equivariant cyclic homology. Assume that $ V $ is an essential $ H $-module and that $ A $ is an $ H $-algebra. 
Moreover let $ b: V \times V \rightarrow \mathbb{C} $ be an equivariant bounded linear map. We define 
an $ H $-algebra $ l(b;A) $ by equipping the space $ V \cotimes A \cotimes V $ with the multiplication 
$$
(v_1 \otimes a_1 \otimes w_1)(v_2 \otimes a_2 \otimes w_2) 
= b(w_1, v_2)\, v_1 \otimes a_1a_2 \otimes w_2 
$$ 
and the diagonal $ H $-action. \\
As a particular case of this construction consider the space $ V = \hat{H} $ with the regular action of $ H $ given by 
$ (t \hit f)(x) = f(xt) $ and the pairing 
$$
\beta(f,g) = \hat{\psi}(fg). 
$$
We write $ \mathcal{K}_H $ for the algebra $ l(\beta; \mathbb{C}) $ and 
$ A \cotimes \mathcal{K}_H $ for $ l(\beta; A) $. Remark that the action on $ A \cotimes \mathcal{K}_H $ is not the diagonal 
action in general. We denote an element $ f \otimes a \otimes g $ in this algebra
by $ |f \ket \otimes a  \otimes \bra g| $ in the sequel. 
Using the isomorphism $ \hat{\mathcal{F}}_r S^{-1}: \hat{H} \rightarrow H $ we identify the above 
pairing with a pairing $ H \times H \rightarrow \mathbb{C} $. The corresponding action of $ H $ on itself is given 
by left multiplication and using the normalization $ \phi = S(\psi) $ we obtain the formula
\begin{align*}
\beta(x,y) = \beta(S \mathcal{G}_r S(x), S \mathcal{G}_r S(y)) = \beta(\mathcal{F}_l(x), \mathcal{F}_l(y)) = \phi(S^{-1}(y)x) = \psi(S(x) y)
\end{align*}
for the above pairing expressed in terms of $ H $. \\
Let $ H $ be a bornological quantum group and let $ A $ be an $ H $-algebra. 
We define a bounded linear map $ \gamma_A: A \rtimes H \rtimes \hat{H} \rightarrow A \cotimes \mathcal{K}_H $ by 
\begin{equation*}
\gamma_A(a \rtimes x \rtimes \mathcal{F}_l(y)) = |y_{(1)} S(x_{(2)}) \ket \otimes y_{(2)} S(x_{(1)}) \cdot a \otimes \bra y_{(3)}|
\end{equation*}
and it is straightforward to check that $ \gamma_A $ is a bornological isomorphism. Using 
$$
\mathcal{F}_l(y^1)_{(1)} \hit x^2 \otimes \mathcal{F}_l(y^1)_{(2)} \mathcal{F}_l(y^2) = 
x^2_{(1)} \otimes \phi(x^2_{(2)} S^{-1}(y^2_{(1)})y^1) \mathcal{F}_l(y^2_{(2)}) 
$$
we compute 
\begin{align*}
\gamma_A((a^1 \rtimes &x^1 \rtimes \mathcal{F}_l(y^1))(a^2 \rtimes x^2 \rtimes \mathcal{F}_l(y^2))) \\
&= \gamma_A((a^1 \rtimes x^1)(a^2 \rtimes x^2_{(1)}) \phi(x^2_{(2)} S^{-1}(y^2_{(1)}) y^1) \rtimes \mathcal{F}_l(y^2_{(2)})) \\
&= \gamma_A(a^1 x^1_{(1)} \cdot a^2 \rtimes x^1_{(2)} x^2_{(1)} \phi(x^2_{(2)} S^{-1}(y^2_{(1)}) y^1) \rtimes \mathcal{F}_l(y^2_{(2)})) \\
&= |y^2_{(2)} S(x^1_{(3)} x^2_{(2)}) \ket \otimes y^2_{(3)}S(x^1_{(2)} x^2_{(1)}) \cdot (a^1  x^1_{(1)} \cdot a^2) \\ 
& \qquad \otimes \phi(x^2_{(3)} S^{-1}(y^2_{(1)}) y^1) \bra y^2_{(4)}| \\
&= |y^2_{(2)} S(x^1_{(2)} x^2_{(3)}) \ket \otimes (y^2_{(3)}S(x^1_{(1)} x^2_{(2)}) \cdot a^1) (y^2_{(4)} S(x^2_{(1)}) \cdot a^2) \\
&\qquad \otimes \phi(x^2_{(4)} S^{-1}(y^2_{(1)}) y^1) \bra y^2_{(5)}| 
\end{align*}
and since $ \phi $ is left invariant this is equal to
\begin{align*}
&= |y^2_{(4)} S(x^1_{(2)} S^{-1}(S^{-1}(y^2_{(3)}) y^1_{(1)})) \ket \otimes \\
&\qquad (y^2_{(5)} S(x^1_{(1)} S^{-1}(S^{-1}(y^2_{(2)}) y^1_{(2)})) \cdot a^1)(y^2_{(6)} S(x^2_{(1)}) \cdot a^2) \\
&\qquad \otimes \phi(x^2_{(2)} S^{-1}(y^2_{(1)}) y^1_{(3)}) \bra y^2_{(7)}| \\
&= |y^1_{(1)} S(x^1_{(2)}) \ket \otimes (y^1_{(2)} S(x^1_{(1)}) \cdot a^1)(y^2_{(2)} S(x^2_{(1)}) \cdot a^2) \\
&\qquad \otimes \phi(x^2_{(2)} S^{-1}(y^2_{(1)})y^1_{(3)}) \bra y^2_{(3)}| \\
&= |y^1_{(1)} S(x^1_{(2)}) \ket \otimes (y^1_{(2)} S(x^1_{(1)}) \cdot a^1)(y^2_{(2)} S(x^2_{(1)}) \cdot a^2) \\
&\qquad \otimes \psi(S(y^1_{(3)}) y^2_{(1)} S(x^2_{(2)})) \bra y^2_{(3)}| \\
&= \gamma_A(a^1 \rtimes x^1 \rtimes \mathcal{F}_l(y^1)) \gamma_A(a^2 \rtimes x^2 \rtimes \mathcal{F}_l(y^2))
\end{align*}
where we use $ \phi = S(\psi) $. 
It follows that $ \gamma_A $ is an algebra homomorphism. In addition it is easily seen that $ \gamma_A $ is equivariant. 
Consequently we obtain the following analogue of the Takesaki-Takai duality theorem. 
\begin{prop}\label{TakTak}
Let $ H $ be a bornological quantum group and let $ A $ be an $ H $-algebra. Then the map $ \gamma_A: A \rtimes H \rtimes \hat{H} 
\rightarrow A \cotimes \mathcal{K}_H $ is an equivariant algebra isomorphism. 
\end{prop} 
For algebraic quantum groups a discussion of Takesaki-Takai duality is contained in \cite{DvDZ}. 
More information on similar duality results in the context of Hopf algebras can be found in \cite{Montgomery}. \\
If $ H = \D(G) $ is the smooth convolution algebra of a locally compact group $ G $ then an $ H $-algebra is 
the same thing as a $ G $-algebra. As a special case of proposition \ref{TakTak} one obtains that 
for every $ G $-algebra $ A $ the double crossed product $ A \rtimes H \rtimes \hat{H} $ is isomorphic to the $ G $-algebra 
$ A \cotimes \mathcal{K}_G $ used in \cite{Voigtepch}. 

\section{Radford's formula} \label{secradford}

In this section we prove a formula for the fourth power of the antipode in terms of the modular elements 
of a bornological quantum group and its dual. This formula was obtained by Radford in the setting of finite 
dimensional Hopf algebras \cite{Radford}.\\
Let $ H $ be a bornological quantum group. If $ \phi $ is a left Haar functional on $ H $ 
there exists a unique multiplier $ \delta \in M(H) $ such that 
$$
(\phi \cotimes \id)\Delta(x) = \phi(x) \delta
$$ 
for all $ x \in H $. The multiplier $ \delta $ is called the modular element of $ H $ and measures the failure of $ \phi $ 
from being right invariant. It is shown in \cite{Voigtbqg} that $ \delta $ is invertible 
with inverse $ S(\delta) = S^{-1}(\delta) = \delta^{-1} $ and that one has 
$ \Delta(\delta)= \delta \otimes \delta $ as well as $ \epsilon(\delta) = 1 $. In terms of the 
dual quantum group the modular element $ \delta $ defines a character, that is, an essential homomorphism 
from $ \hat{H} $ to $ \mathbb{C} $. Similarly, there exists a unique modular element $ \hat{\delta} \in M(\hat{H}) $ for the dual quantum group 
which satisfies 
$$ 
(\hat{\phi} \cotimes \id)\hat{\Delta}(f) = \hat{\phi}(f) \hat{\delta} 
$$ 
for all $ f \in \hat{H} $. \\
The Haar functionals of a bornological quantum group are uniquely determined up to a scalar multiple. In many situations it is convenient to 
fix a normalization at some point. However, in the discussion below it is not necessary to keep track of the scaling of the 
Haar functionals. If $ \omega $ and $ \eta $ are linear functionals we shall write $ \omega \equiv \eta $ if there exists 
a nonzero scalar $ \lambda $ such that $ \omega = \lambda \eta $. We use the same notation for elements in a bornological 
quantum group or linear maps that differ by some nonzero scalar multiple. Moreover we shall identify $ H $ with its double dual 
using Pontrjagin duality. \\
To begin with observe that the bounded linear functional $ \delta \hit \phi $ on $ H $ defined by 
$$
(\delta \hit \phi)(x) = \phi(x \delta) 
$$
is faithful and satisfies 
\begin{align*}
((\delta \hit \phi) \cotimes \id)&\Delta(x) = (\phi \cotimes \id)(\Delta(x \delta) (1 \otimes \delta^{-1})) = 
\phi(x \delta) \delta \delta^{-1} = (\delta \hit \phi)(x). 
\end{align*}
It follows that  $ \delta \hit \phi $ is a right Haar functional on $ H $. In a similar way we obtain a right Haar functional 
$ \phi \hitby \delta $ on $ H $. Hence 
$$ 
\delta \hit \phi \equiv \psi \equiv \phi \hitby \delta
$$ 
by uniqueness of the right Haar functional which yields in particular the relations
$$
\mathcal{F}_l(x\delta) \equiv\mathcal{G}_l(x), \qquad \mathcal{F}_r(\delta x) \equiv \mathcal{G}_r(x)
$$
for the Fourier transform. \\
According to Pontrjagin duality we have $ x = \hat{\mathcal{G}}_l \mathcal{F}_l S(x) $ for all $ x \in H $ and using 
$ \hat{\mathcal{G}}_l(f) \equiv \hat{\mathcal{F}}_l(f \hat{\delta}) $ for every $ f \in \hat{H} $ as well as 
\begin{align*}
(\mathcal{F}_l(S(&x))\hat{\delta})(h) = \phi(h_{(1)} S(x)) \hat{\delta}(h_{(2)}) = \phi(h S(x_{(2)})) \hat{\delta}(\delta S^2(x_{(1)})) \\
&\equiv \mathcal{F}_l(S(x_{(2)}))(h) \hat{\delta}(x_{(1)}) = \mathcal{F}_l(S(x \hitby \hat{\delta}))(h)
\end{align*}
we obtain $ x \equiv \hat{\mathcal{F}}_l \mathcal{F}_l S(x \hitby \hat{\delta}) $ or equivalently 
\begin{equation} \label{radford1h1}
S^{-1}(\hat{\delta} \hit x) \equiv \hat{\mathcal{F}}_l \mathcal{F}_l(x). 
\end{equation}
Using equation (\ref{radford1h1}) and the formula $ x \equiv S^{-1} \hat{\mathcal{F}}_l \mathcal{F}_r(x) $ obtained from Pontrjagin duality 
we compute 
$$
\hat{\mathcal{F}}_l \mathcal{F}_l(\hat{\delta}^{-1} \hit S^2(x)) \equiv S^{-1} S^2(x) = S(x) \equiv \hat{\mathcal{F}}_l \mathcal{F}_r(x) 
$$
which implies 
\begin{equation} \label{radford1}
\mathcal{F}_l(S^2(x)) \equiv \mathcal{F}_r(\hat{\delta} \hit x)
\end{equation}
since $ \hat{\mathcal{F}}_l $ is an isomorphism. 
Similarly, we have $ \hat{\mathcal{F}}_l S \mathcal{G}_l \equiv \id $ and using $ \hat{\mathcal{F}}_l(f) \equiv \hat{\mathcal{G}}_l(f \hat{\delta}^{-1}) $ 
together with 
\begin{align*}
(S\mathcal{G}_l(&x) \hat{\delta}^{-1})(h) = \psi(S(h_{(1)}) x) \hat{\delta}^{-1}(h_{(2)}) \\
&\equiv \psi(S(h) x_{(2)}) \hat{\delta}^{-1}(x_{(1)}) = S\mathcal{G}_l(x \hitby \hat{\delta}^{-1})(h)
\end{align*}
we obtain $ \hat{\mathcal{G}}_l S \mathcal{G}_l(x) \equiv x \hitby \hat{\delta} $.  
According to the relation $ \mathcal{F}_l(x \delta) \equiv \mathcal{G}_l(x) $ this may be rewritten as 
$ S^{-1}\hat{\mathcal{F}}_r \mathcal{F}_l(x \delta) \equiv x \hitby \hat{\delta} $ 
which in turn yields 
\begin{equation} \label{radford2h1}
\hat{\mathcal{F}}_r \mathcal{F}_l(x) \equiv S((x \delta^{-1}) \hitby \hat{\delta}). 
\end{equation}
Due to Pontrjagin duality we have $ x = \hat{\mathcal{F}}_r \mathcal{G}_r S(x) $ for all $ x \in H $ and 
using $ \mathcal{G}_r(x) \equiv \mathcal{F}_r(\delta x) $ we obtain 
\begin{equation} \label{radford2h2}
\hat{\mathcal{F}}_r \mathcal{F}_r(S(x)) \equiv \hat{\mathcal{F}}_r \mathcal{G}_r(\delta^{-1} S(x)) = x \delta. 
\end{equation}
According to equation (\ref{radford2h1}) and equation (\ref{radford2h2}) we have
\begin{align*}
\hat{\mathcal{F}}_r& \mathcal{F}_l(S^{-2}(\delta^{-1}(x\hitby \hat{\delta}^{-1}) \delta)) \equiv 
S((S^{-2}(\delta^{-1}(x \hitby \hat{\delta}^{-1}) \delta)\delta^{-1}) \hitby \hat{\delta}) \\
&\equiv S(S^{-2}(\delta^{-1}x)) = S^{-1}(\delta^{-1}x) \equiv \hat{\mathcal{F}}_r \mathcal{F}_r(x)
\end{align*}
and since $ \hat{\mathcal{F}}_r $ is an isomorphism this implies
\begin{equation} \label{radford2}
\mathcal{F}_r(S^2(x)) \equiv \mathcal{F}_l(\delta^{-1}(x \hitby \hat{\delta}^{-1}) \delta).
\end{equation}
Assembling these relations we obtain the following result. 
\begin{prop}\label{radford}
Let $ H $ be a bornological quantum group and let $ \delta $ and $ \hat{\delta} $ be the modular 
elements of $ H $ and $ \hat{H} $, respectively. Then 
$$
S^4(x) = \delta^{-1} (\hat{\delta} \hit x \hitby \hat{\delta}^{-1}) \delta
$$
for all $ x \in H $. 
\end{prop}
\proof Using equation (\ref{radford1}) and equation (\ref{radford2}) we compute 
$$
\mathcal{F}_l(S^4(x)) \equiv \mathcal{F}_r(\hat{\delta} \hit S^2(x)) = 
\mathcal{F}_r(S^2(\hat{\delta} \hit x)) \equiv \mathcal{F}_l(\delta^{-1}(\hat{\delta} \hit x \hitby \hat{\delta}^{-1}) \delta)
$$
which implies 
$$
S^4(x) \equiv \delta^{-1}(\hat{\delta} \hit x \hitby \hat{\delta}^{-1}) \delta
$$
for all $ x \in H $ since $ \mathcal{F}_l $ is an isomorphism. The claim follows from the observation that both sides of the previous 
equation define algebra automorphisms of $ H $. \qed 

\section{Anti-Yetter-Drinfeld modules} \label{seccovmod}

In this section we introduce the notion of an anti-Yetter-Drinfeld module over a bornological 
quantum group. Moreover we discuss the concept of a paracomplex. \\
We begin with the definition of an anti-Yetter-Drinfeld module. In the context of Hopf algebras this 
notion was introduced in \cite{HKRS}. 
\begin{definition}
Let $ H $ be a bornological quantum group. An $ H $-anti-Yetter-Drinfeld module is an essential left $ H $-module $ M $ 
which is also an essential left $ \hat{H} $-module such that 
$$
t \cdot (f \cdot m) = (S^2(t_{(1)}) \hit f \hitby S^{-1}(t_{(3)})) \cdot (t_{(2)} \cdot m).
$$
for all $ t \in H, f \in \hat{H} $ and $ m \in M $. 
A homomorphism $ \xi: M \rightarrow N $ between anti-Yetter-Drinfeld modules 
is a bounded linear map which is both $ H $-linear and $ \hat{H} $-linear. 
\end{definition}
We will not always mention explicitly the underlying bornological quantum group when 
dealing with anti-Yetter-Drinfeld modules. Moreover we shall use the abbreviations $ \ayd $-module and $ \ayd $-map for 
anti-Yetter-Drinfeld modules and their homomorphisms. \\
According to theorem \ref{modcomod} a left $ \hat{H} $-module structure corresponds to a right $ H $-comodule structure. 
Hence an $ \ayd $-module can be described equivalently as a bornological vector space $ M $ equipped with an 
essential $ H $-module structure and an $ H $-comodule structure satisfying a certain compatibility condition. 
Formally, this compatibility condition can be written down as 
\begin{equation*}
(t \cdot m)_{(0)} \otimes (t \cdot m)_{(1)} = 
t_{(2)} \cdot m_{(0)} \otimes t_{(3)} m_{(1)} S(t_{(1)})
\end{equation*}
for all $ t \in H $ and $ m \in M  $. \\
We want to show that $ \ayd $-modules can be interpreted as essential modules over a 
certain bornological algebra. Following the notation in \cite{HKRS} this algebra will be denoted by $ \AYD(H) $. 
As a bornological vector space we have $ \AYD(H) = \hat{H} \cotimes H $ and the multiplication is defined by the formula
$$
(f \otimes x)\cdot(g \otimes y) = f (S^2(x_{(1)}) \hit g \hitby S^{-1}(x_{(3)})) \otimes x_{(2)} y. 
$$
There exists an algebra homomorphism $ \iota_H: H \rightarrow M(\AYD(H)) $ given by 
$$
\iota_H(x) \cdot (g \otimes y) = S^2(x_{(1)}) \hit g \hitby S^{-1}(x_{(3)}) \otimes x_{(2)} y
$$
and 
$$
(g \otimes y) \cdot \iota_H(x) = g \otimes yx. 
$$
It is easily seen that $ \iota_H $ is injective. 
Similarly, there is an injective algebra homomorphism $ \iota_{\hat{H}}: \hat{H} \rightarrow M(\AYD(H)) $ 
given by 
$$
\iota_{\hat{H}}(f) \cdot (g \otimes y) = fg \otimes y
$$
as well as 
$$
(g \otimes y) \cdot \iota_{\hat{H}}(f) = g(S^2(y_{(1)}) \hit f \hitby S^{-1}(y_{(3)})) \otimes y_{(2)}
$$
and we have the following result. 
\begin{prop} \label{AYDHess}
For every bornological quantum group $ H $ the bornological algebra $ \AYD(H) $ is essential. 
\end{prop} 
\proof The homomorphism $ \iota_{\hat{H}} $ induces on $ \AYD(H) $ the structure of an essential left 
$ \hat{H} $-module. Similarly, the space $ \AYD(H) $ becomes an essential left $ H $-module 
using the homomorphism $ \iota_H $. In fact, if we write $ \hat{H}_\tau $ for the space $ \hat{H} $ equipped with the 
trivial $ H $-action the map $ c: \AYD(H) \rightarrow \hat{H}_\tau \cotimes H $ given by 
$$ 
c(f \otimes x) = S(x_{(1)}) \hit f \hitby x_{(3)} \otimes x_{(2)}
$$
is an $ H $-linear isomorphism. Actually, the actions of $ \hat{H} $ and $ H $ on $ \AYD(H) $ are defined in such a way that 
$ \AYD(H) $ becomes an $ \ayd $-module. 
Using the canonical isomorphism $ H \cotimes_H \AYD(H) \cong \AYD(H) $ one obtains an essential $ \hat{H} $-module structure 
on $ H \cotimes_H \AYD(H) $ given explicitly by the formula 
$$
f \cdot (x \otimes g \otimes y) = x_{(2)} \otimes (S(x_{(1)}) \hit f \hitby x_{(3)}) g \otimes y.  
$$
Correspondingly, we obtain a natural isomorphism $ \hat{H} \cotimes_{\hat{H}} (H \cotimes_H \AYD(H)) \cong \AYD(H) $. 
It is straighforward to verify that the identity map $ \hat{H} \cotimes (H \cotimes \AYD(H)) \cong 
\hat{H} \cotimes H \cotimes \hat{H} \cotimes H \rightarrow \AYD(H) \cotimes \AYD(H) $ induces an 
isomorphism $ \hat{H} \cotimes_{\hat{H}} (H \cotimes_H \AYD(H)) 
\rightarrow \AYD(H) \cotimes_{\AYD(H)} \AYD(H) $. This yields the claim.  \qed \\
We shall now characterize $ \ayd $-modules as essential modules over $ \AYD(H) $. 
\begin{prop}\label{covrepprop} Let $ H $ be a bornological quantum group. Then the 
category of $ \ayd $-modules over $ H $ is isomorphic to the category of essential left $ \AYD(H) $-modules. 
\end{prop}
\proof Let $ M \cong \AYD(H) \cotimes_{\AYD(H)} M $ be an essential $ \AYD(H) $-module. 
Then we obtain a left $ H $-module structure and a left $ \hat{H} $-module structure on $ M $ 
using the canonical homomorphisms $ \iota_H: H \rightarrow M(\AYD(H)) $ and 
$ \iota_{\hat{H}}: \hat{H} \rightarrow M(\AYD(H)) $. 
Since the action of $ H $ on $ \AYD(H) $ is essential we have natural isomorphisms
$$ 
H \cotimes_H M \cong H \cotimes_H \AYD(H) \cotimes_{\AYD(H)} M \cong 
\AYD(H) \cotimes_{\AYD(H)} M \cong M 
$$ 
and hence $ M $ is an essential $ H $-module. Similarly we have 
$$
\hat{H} \cotimes_{\hat{H}} M \cong  \hat{H} \cotimes_{\hat{H}} \AYD(H) \cotimes_{\AYD(H)} \cotimes M \cong 
\AYD(H) \cotimes_{\AYD(H)} M \cong M 
$$
since $ \AYD(H) $ is an essential $ \hat{H} $-module. These module actions yield the structure of 
an $ \ayd $-module on $ M $. \\
Conversely, assume that $ M $ is an $ H $-$ \ayd $-module. Then we obtain 
an $ \AYD(H) $-module structure on $ M $ by setting 
$$
(f \otimes t) \cdot m = f \cdot (t \cdot m) 
$$
for $ f \in \hat{H} $ and $ t \in H $. Since $ M $ is an essential $ H $-module 
we have a natural isomorphism $ H \cotimes_H M \cong M $. As in the proof of 
proposition \ref{AYDHess} we obtain an induced essential $ \hat{H} $-module structure on 
$ H \cotimes_H M $ and canonical isomorphisms $ \AYD(H) \cotimes_{\AYD(H)} M \cong \hat{H} \cotimes_{\hat{H}} (H \cotimes_H M) \cong M $. 
It follows that $ M $ is an essential $ \AYD(H) $-module. \\
The previous constructions are compatible with morphisms and it is easy to check that they are 
inverse to each other. This yields the assertion. \qed \\
There is a canonical operator $ T $ on every $ \ayd $-module which plays a crucial role in equivariant 
cyclic homology. In order to define this operator it is convenient to pass from $ \hat{H} $ to $ H $ in the 
first tensor factor of $ \AYD(H) $. More precisely, consider the bornological isomorphism 
$ \lambda: \AYD(H) \rightarrow H \cotimes H $ given by 
$$
\lambda(f \otimes y) = \hat{\mathcal{F}}_l(f) \otimes y \hitby \hat{\delta}^{-1}
$$ 
where $ \hat{\delta} \in M(\hat{H}) $ is the modular function of $ \hat{H} $. The inverse map is given by 
$$
\lambda^{-1}(x \otimes y) = S \mathcal{G}_l(x) \otimes y \hitby \hat{\delta}. 
$$
It is straightforward to check that the left $ H $-action on $ \AYD(H) $ corresponds to 
$$
t \cdot (x \otimes y) = t_{(3)} x S(t_{(1)}) \otimes t_{(2)} y
$$
and the left $ \hat{H} $-action becomes 
$$
f \cdot (x \otimes y) = (f \hit x) \otimes y
$$
under this isomorphism. The right $ H $-action is identified with 
$$
(x \otimes y) \cdot t = x \otimes y (t \hitby \hat{\delta}^{-1})
$$
and the right $ \hat{H} $-action corresponds to 
$$
(x \otimes y) \cdot g = x_{(2)} (S^2(y_{(2)}) \hit g \hitby S^{-1}(y_{(4)}))(S^{-2}(x_{(1)}) \delta) \otimes \hat{\delta}(y_{(1)}) y_{(3)} \hitby \hat{\delta}^{-1}
$$
where $ \delta $ is the modular function of $ H $. 
Using this description of $ \AYD(H) $ we obtain the following result. 
\begin{prop} \label{Toplemma} The bounded linear map $ T: \AYD(H) \rightarrow \AYD(H) $ defined by 
$$
T(x \otimes y) = x_{(2)} \otimes S^{-1}(x_{(1)}) y
$$
is an isomorphism of $ \AYD(H) $-bimodules. 
\end{prop}
\proof It is evident that $ T $ is a bornological isomorphism with inverse given by 
$ T^{-1}(x \otimes y) = x_{(2)} \otimes x_{(1)}y $. We compute 
\begin{align*}
T(t \cdot (x &\otimes y)) = T(t_{(3)} x S(t_{(1)}) \otimes t_{(2)} y) \\
&= t_{(5)} x_{(2)} S(t_{(1)}) \otimes S^{-1}(t_{(4)} x_{(1)} S(t_{(2)}))t_{(3)} y \\
&= t_{(3)} x_{(2)} S(t_{(1)}) \otimes t_{(2)} S^{-1}(x_{(1)}) y \\
&= t \cdot T(x \otimes y)
\end{align*}
and it is clear that $ T $ is left $ \hat{H} $-linear. Consequently $ T $ is a left $ \AYD(H) $-linear map. 
Similarly, we have 
$$
T((x \otimes y)\cdot t) = T(x \otimes y(t \hitby \hat{\delta}^{-1})) = x_{(2)} \otimes S^{-1}(x_{(1)}) y(t \hitby \hat{\delta}^{-1}) 
= T(x \otimes y) \cdot t
$$
and hence $ T $ is right $ H $-linear. In order to prove that $ T $ is right $ \hat{H} $-linear we compute 
\begin{align*}
T&^{-1}((x \otimes y)\cdot g) \\
&= T^{-1}(x_{(2)} (S^2(y_{(2)}) \hit g \hitby S^{-1}(y_{(4)}))(S^{-2}(x_{(1)}) \delta) \otimes \hat{\delta}(y_{(1)}) y_{(3)} \hitby \hat{\delta}^{-1}) \\
&= x_{(3)} (S^2(y_{(2)}) \hit g \hitby S^{-1}(y_{(4)}))(S^{-2}(x_{(1)}) \delta) \otimes \hat{\delta}(y_{(1)}) x_{(2)} (y_{(3)} \hitby \hat{\delta}^{-1}) \\
&= x_{(3)} (S^2(y_{(2)}) \hit g \hitby S^{-1}(y_{(4)}))(S^{-2}(\hat{\delta} \hit x_{(1)}) \delta) \otimes 
\hat{\delta}(y_{(1)}) (x_{(2)} y_{(3)}) \hitby \hat{\delta}^{-1} 
\end{align*}
which according to proposition \ref{radford} is equal to
\begin{align*}
&= x_{(3)} (S^2(y_{(2)}) \hit g \hitby S^{-1}(y_{(4)}))(\delta S^2(x_{(1)} \hitby \hat{\delta})) \otimes 
\hat{\delta}(y_{(1)}) (x_{(2)} y_{(3)}) \hitby \hat{\delta}^{-1}  \\
&= x_{(6)} (S^2(x_{(2)}) S^2(y_{(2)}) \hit g \hitby S^{-1}(y_{(4)}) S^{-1}(x_{(4)}))(S^{-2}(x_{(5)}) \delta) \\
&\qquad \otimes \hat{\delta}(x_{(1)} y_{(1)}) (x_{(3)} y_{(3)}) \hitby \hat{\delta}^{-1} \\
&= (x_{(2)} \otimes x_{(1)} y) \cdot g \\
&= T^{-1}(x \otimes y) \cdot g.
\end{align*}
We conclude that $ T $ is a right $ \AYD(H) $-linear map as well. \qed \\ 
If $ M $ is an arbitrary $ \ayd $-module we define $ T: M \rightarrow M $ by 
$$
T(F \otimes m) = T(F) \otimes m 
$$
for $ F \otimes m \in \AYD(H) \cotimes_{\AYD(H)} M $. Due to proposition \ref{covrepprop} and lemma \ref{Toplemma} this definition makes 
sense. 
\begin{prop} \label{covparaadd} 
The operator $ T $ defines a natural isomorphism $ T: \id \rightarrow \id $ of the identity functor 
on the category of $ \ayd $-modules. 
\end{prop}
\proof It is clear from the construction that $ T: M \rightarrow M $ is an isomorphism for all $ M $. 
If $ \xi: M \rightarrow N $ is an $ \ayd $-map the equation $ T \xi = \xi T $ follows easily 
after identifying $ \xi $ with the map $ \id \cotimes \xi: \AYD(H) \cotimes_{\AYD(H)} M 
\rightarrow \AYD(H) \cotimes_{\AYD(H)} N $. This yields the assertion. \qed \\
If the bornological quantum group $ H $ is unital one may construct the operator $ T $ on an $ \ayd $-module $ M $ directly from the 
coaction $ M \rightarrow M \cotimes H $ corresponding to the action of $ \hat{H} $. More precisely, one has the formula
$$
T(m) = S^{-1}(m_{(1)}) \cdot m_{(0)}
$$
for every $ m \in M $. \\
Using the terminology of \cite{Voigtepch} it follows from proposition \ref{covparaadd} that the category of $ \ayd $-modules 
is a para-additive category in a natural way. This leads in particular to the concept of a paracomplex of $ \ayd $-modules. 
\begin{definition} 
A paracomplex $ C = C_0 \oplus C_1 $ consists of $ \ayd $-modules $ C_0 $ and $ C_1 $ 
and $ \ayd $-maps $ \partial_0: C_0 \rightarrow C_1 $ and $ \partial_1: C_1 \rightarrow C_0 $ such that
$$ 
\partial^2 = \id - T 
$$
where the differential $ \partial: C \rightarrow C_1 \oplus C_0 \cong C $ is the composition of $ \partial_0 \oplus \partial_1 $
with the canonical flip map. 
\end{definition}
The morphism $ \partial $ in a paracomplex is called a differential although it usually 
does not satisfy the relation $ \partial^2 = 0 $. As for ordinary complexes one defines chain maps between paracomplexes and homotopy equivalences.  
We always assume that such maps are compatible with the $ \ayd $-module structures. Let us point out that it does not make sense to 
speak about the homology of a paracomplex in general. \\
The paracomplexes we will work with arise from paramixed complexes in the following sense.
\begin{definition} 
A paramixed complex $ M $ is a sequence of
$ \ayd $-modules $ M_n $ together with $ \ayd $-maps $ b $ of degree $ -1 $ and $ B $ of degree $ + 1 $ satisfying
$ b^2 = 0 $, $ B^2 = 0 $ and
\begin{equation*}
[b,B] = bB + Bb = \id - T.
\end{equation*}
\end{definition}
If $ T $ is equal to the identity operator on $ M $ this reduces of course to the definition of a mixed complex. 

\section{Equivariant differential forms}

In this section we define equivariant differential forms and the equivariant $ X $-complex. 
Moreover we discuss the properties of the periodic tensor algebra of an $ H $-algebra. 
These are the main ingredients in the construction of equivariant cyclic homology. \\
Let $ H $ be a bornological quantum group. If $ A $ is an $ H $-algebra we obtain a left 
action of $ H $ on the space $ H \cotimes \Omega^n(A) $ by 
\begin{equation*}
t \cdot (x \otimes \omega) = t_{(3)} x S(t_{(1)}) \otimes t_{(2)} \cdot \omega
\end{equation*}
for $ t, x \in H $ and $ \omega \in \Omega^n(A) $. Here $ \Omega^n(A) = A^+ \cotimes A^{\cotimes n} $ for $ n > 0 $ is the space 
of noncommutative $ n $-forms over $ A $ with the diagonal $ H $-action. For $ n = 0 $ one defines $ \Omega^0(A) = A $. 
There is a left action of the dual quantum group $ \hat{H} $ on $ H \cotimes \Omega^n(A) $ given by 
\begin{equation*}
f\cdot (x \otimes \omega) = (f \hit x) \otimes \omega = f(x_{(2)}) x_{(1)} \otimes \omega.
\end{equation*}
By definition, the equivariant $ n $-forms $ \Omega^n_H(A) $ are the space $ H \cotimes \Omega^n(A) $ 
together with the $ H $-action and the $ \hat{H} $-action described above. We compute 
\begin{align*}
t \cdot (f \cdot &(x \otimes \omega)) = t \cdot (f(x_{(2)}) x_{(1)} \otimes \omega) \\
&= f(x_{(2)}) t_{(3)} x_{(1)} S(t_{(1)}) \otimes t_{(2)} \cdot \omega \\
&= (S^2(t_{(1)}) \hit f \hitby S^{-1}(t_{(5)})) \cdot (t_{(4)} x S(t_{(2)}) \otimes t_{(3)} \cdot \omega) \\
&= (S^2(t_{(1)}) \hit f \hitby S^{-1}(t_{(3)})) \cdot (t_{(2)} \cdot (x \otimes \omega))
\end{align*}
and deduce that $ \Omega^n_H(A) $ is an $ H $-$ \ayd $-module. We let $ \Omega_H(A) $ be the direct 
sum of the spaces $ \Omega^n_H(A) $. \\
Now we define operators $ d $ and $ b_H $ on $ \Omega_H(A) $ by
\begin{equation*}
d(x \otimes \omega) = x \otimes d\omega
\end{equation*}
and
\begin{equation*}
b_H(x \otimes \omega da) = (-1)^{|\omega|} (x \otimes \omega a -
x_{(2)} \otimes (S^{-1}(x_{(1)}) \cdot a) \omega).
\end{equation*}
The map $ b_H $ should be thought of as a twisted version of the usual
Hochschild operator. We compute 
\begin{align*}
b_H^2&(x \otimes \omega da db) = (-1)^{|\omega| + 1} b_H(x \otimes \omega da b - 
x_{(2)} \otimes (S^{-1}(x_{(1)}) \cdot b) \omega da) \\
&= (-1)^{|\omega| + 1} b_H(x \otimes \omega d(ab) 
- x \otimes \omega a db - x_{(2)} \otimes (S^{-1}(x_{(1)}) \cdot b) \omega da) \\
&= - (x \otimes \omega ab - x_{(2)} \otimes S^{-1}(x_{(1)}) \cdot (ab) \omega  
- x \otimes \omega ab + x_{(2)} \otimes (S^{-1}(x_{(1)}) \cdot b) \omega a \\
&\qquad - x_{(2)} \otimes (S^{-1}(x_{(1)}) \cdot b) \omega a 
+ x_{(2)} \otimes S^{-1}(x_{(1)}) \cdot (ab) \omega) = 0
\end{align*}
which shows that $ b_H^2 $ is a differential as in the nonequivariant situation. 
Let us discuss the compatibility of $ d $ and $ b_H $ with the
$ \ayd $-module structure. It is easy to check that $ d $ is an $ \ayd $-map and 
that the operator $ b_H $ is $ \hat{H} $-linear. Moreover we compute
\begin{align*}
&b_H(t \cdot(x \otimes \omega da)) = (-1)^{|\omega|} (t_{(4)}x S(t_{(1)}) \otimes (t_{(2)} \cdot \omega) (t_{(3)} \cdot a) \\
&\qquad - t_{(6)} x_{(2)} S(t_{(1)}) \otimes
(S^{-1}(t_{(5)} x_{(1)} S(t_{(2)})) t_{(4)} \cdot a) (t_{(3)} \cdot \omega)) \\
&= (-1)^{|\omega|} (t_{(3)}x S(t_{(1)}) \otimes t_{(2)} \cdot (\omega a) - t_{(4)} x_{(2)} S(t_{(1)}) \otimes (t_{(2)}
S^{-1}(x_{(1)}) \cdot a) (t_{(3)}\cdot \omega)) \\
&= t \cdot b_H(x \otimes \omega da)
\end{align*}
and deduce that $ b_H $ is an $ \ayd $-map as well. \\
Similar to the non-equivariant case we use $ d $ and $ b_H $ to
define an equivariant Karoubi operator $ \kappa_H $ and an
equivariant Connes operator $ B_H $ by 
\begin{equation*}
\kappa_H = 1 - (b_H d + d b_H)
\end{equation*}
and 
\begin{equation*}
B_H = \sum_{j = 0}^n \kappa_H^j d
\end{equation*}
on $ \Omega^n_H(A) $. Let us record the following explicit formulas. 
For $ n > 0 $ we have
\begin{equation*}
\kappa_H(x \otimes \omega da) = (-1)^{n - 1} x_{(2)} \otimes (S^{-1}(x_{(1)}) \cdot da) \omega
\end{equation*}
on $ \Omega^n_H(A) $ and in addition $ \kappa_H(x \otimes a) = x_{(2)} \otimes S^{-1}(x_{(1)}) \cdot a $ on $ \Omega^0_H(A) $. 
For the Connes operator we compute
\begin{equation*}
B_H(x \otimes a_0da_1 \cdots da_n) = \sum_{i = 0}^n (-1)^{ni} x_{(2)} \otimes S^{-1}(x_{(1)}) \cdot(da_{n + 1 - i} \cdots da_n).
da_0 \cdots da_{n - i}
\end{equation*}
Furthermore, the operator $ T $ is given by
\begin{equation*}
T(x \otimes \omega) = x_{(2)} \otimes S^{-1}(x_{(1)}) \cdot \omega
\end{equation*}
on equivariant differential forms. Observe that all operators constructed so far are $ \ayd $-maps and thus commute with $ T $ according to proposition \ref{covparaadd}. 
\begin{lemma}\label{diffformformel}
On $ \Omega^n_H(A) $ the following relations hold:
\begin{bnum}
\item[a)] $ \kappa_H^{n + 1} d = T d $
\item[b)] $ \kappa^n_H = T + b_H \kappa^n_H d $
\item[c)] $ \kappa^n_H b_H = b_H T $
\item[d)] $ \kappa^{n + 1}_H = (\id - db_H) T $
\item[e)] $ (\kappa^{n + 1}_H - T)(\kappa^n_H - T) = 0 $
\item[f)] $ B_H b_H + b_H B_H = \id - T $
\end{bnum}
\end{lemma}
\proof a) follows from the explicit formula for $ \kappa_H $. b) We compute
\begin{align*}
&\kappa_H^n(x \otimes a_0 da_1 \cdots da_n)
= x_{(2)} \otimes S^{-1}(x_{(1)}) \cdot (da_1 \cdots da_n) a_0 \\
&= x_{(2)} \otimes S^{-1}(x_{(1)}) \cdot (a_0 da_1 \cdots da_n) +
(-1)^n b_H(x_{(2)} \otimes
S^{-1}(x_{(1)}) \cdot (da_1 \cdots da_n) da_0) \\
&= x_{(2)} \otimes S^{-1}(x_{(1)}) \cdot (a_0 da_1 \cdots da_n) +
b_H \kappa^n_H d(x \otimes a_0da_1 \cdots da_n)
\end{align*}
which yields the claim. c) follows by applying $ b_H $ to both
sides of b). d) Apply $ \kappa_H $ to b) and use a). e) is a
consequence of b) and d). f) We compute
\begin{align*}
B_H& b_H + b_H B_H  = \sum_{j = 0}^{n - 1} \kappa_H^j db_H  +
\sum_{j = 0}^n b_H \kappa_H^j d
= \sum_{j = 0}^{n - 1} \kappa_H^j(db_H + b_H d) + \kappa_H^n b_H d \\
&= \id - \kappa_H^n(1 - b_H d) = \id - \kappa_H^n(\kappa_H + db_H) = \id
- T + db_H T - Tdb_H  = \id - T
\end{align*}
where we use d) and b). \qed \\
From the definition of $ B_H $ and the fact that $ d^2 = 0 $ we obtain $ B_H^2 = 0 $. 
Let us summarize this discussion as follows.
\begin{prop} Let $ A $ be an $ H $-algebra. The space $ \Omega_H(A) $
of equivariant differential forms is a paramixed complex in the 
category of $ \ayd $-modules.
\end{prop}
We remark that the definition of $ \Omega_H(A) $ for $ H = \D(G) $ differs slightly from the definition of $ \Omega_G(A) $ 
in \cite{Voigtepch} if the locally compact group $ G $ is not unimodular. However, this does not affect the definition of the 
equivariant homology groups. \\
In the sequel we will drop the subscripts when working with the operators on $ \Omega_H(A) $ introduced 
above. For instance, we shall simply write $ b $ instead of $ b_H $ and $ B $ instead of $ B_H $. \\
The $ n $-th level of the Hodge tower associated to $ \Omega_H(A) $ is defined by 
\begin{equation*}
\theta^n \Omega_H(A) = \bigoplus_{j = 0}^{n - 1} \Omega^j_H(A) \oplus \Omega_H^n(A)/b(\Omega^{n + 1}_H(A)).
\end{equation*}
Using the grading into even and odd forms we see that 
$ \theta^n \Omega_H(A) $ together with the boundary operator $ B + b $ becomes a paracomplex. 
By definition, the Hodge tower $ \theta \Omega_H(A) $ of $ A $ is the projective system 
$ (\theta^n \Omega_H(A))_{n \in \mathbb{N}} $ obtained in this way. \\
From a conceptual point of view it is convenient to work with pro-categories in the sequel.  
The pro-category $ \pro(\mathcal{C}) $ over a category $ \mathcal{C} $ consists of projective systems 
in $ \mathcal{C} $. A pro-$ H $-algebra is simply an algebra in the category $ \pro(H \LSMod) $. For instance, 
every $ H $-algebra becomes a pro-$ H $-algebra by viewing it as a constant projective system. More information 
on the use of pro-categories in connection with cyclic homology can be found in \cite{CQ3}, \cite{Voigtepch}. 
\begin{definition} Let $ A $ be a pro-$ H $-algebra. The equivariant $ X $-complex 
$ X_H(A) $ of $ A $ is the paracomplex $ \theta^1 \Omega_H(A) $. Explicitly, we have  
\begin{equation*}
    \xymatrix{
      {X_H(A) \colon \ }
        {\Omega^0_H(A)\;} \ar@<1ex>@{->}[r]^-{d} &
          {\;\Omega^1_H(A)/ b(\Omega^2_H(A)).} 
            \ar@<1ex>@{->}[l]^-{b} 
               } 
\end{equation*}
\end{definition}
We are interested in particular in the equivariant $ X $-complex of the periodic tensor algebra $ \mathcal{T}A $ of an 
$ H $-algebra $ A $. The periodic tensor algebra $ \mathcal{T}A $ is the even 
part of $ \theta \Omega(A) $ equipped with the Fedosov product given by
\begin{equation*}
\omega \circ \eta = \omega \eta - (-1)^{|\omega|} d\omega d\eta
\end{equation*}
for homogenous forms $ \omega $ and $ \eta $. 
The natural projection $ \theta \Omega(A) \rightarrow A $ 
restricts to an equivariant homomorphism $ \tau_A: \mathcal{T}A \rightarrow A $
and we obtain an extension 
\begin{equation*} 
 \xymatrix{
     \mathcal{J}A \;\; \ar@{>->}[r] &
         \mathcal{T}A \ar@{->>}[r]^{\tau_A} &
           A   
     }
\end{equation*}
of pro-$ H $-algebras. \\
The main properties of the pro-algebras $ \mathcal{T}A $ and $ \mathcal{J}A $ are explained in \cite{Meyerthesis}, \cite{Voigtepch}. 
Let us recall some terminology. We write $ \mu^n: N^{\cotimes n} \rightarrow  N $ for the iterated multiplication in a
pro-$ H $-algebra $ N $. Then $ N $ is called locally nilpotent if for every  
equivariant pro-linear map $ f: N \rightarrow C $ with constant range $ C $ there 
exists $ n \in \mathbb{N} $ such that $ f \mu^n = 0 $. It is straightforward to check that the pro-$ H $-algebra 
$ \mathcal{J}A $ is locally nilpotent. \\
An equivariant pro-linear map $ l: A \rightarrow B $ between pro-$ H $-algebras is called a lonilcur if its curvature 
$ \omega_l: A \cotimes A \rightarrow B $ defined by $ \omega_l(a,b) = l(ab) - l(a)l(b) $ is 
locally nilpotent, that is, if for every equivariant pro-linear map 
$ f: B \rightarrow C $ with 
constant range $ C $ there exists $ n \in \mathbb{N} $ such that 
$ f \mu^n_B \omega^{\cotimes n}_l = 0 $. The term lonilcur is an abbreviation for 
"equivariant pro-linear map with \emph{lo}cally \emph{nil}potent \emph{cur}vature". 
Since $ \mathcal{J}A $ is locally nilpotent the natural map $ \sigma_A: A \rightarrow \mathcal{T}A $ 
is a lonilcur. 
\begin{prop}\label{PeriodicTensorAlg} Let $ A $ be an $ H $-algebra. The 
pro-$ H $-algebra $ \mathcal{T}A $ and the lonilcur $ \sigma_A: A 
\rightarrow \mathcal{T}A $ satisfy the following universal property. 
If $ l: A \rightarrow B $ is a lonilcur into a pro-$ H $-algebra $ B $ 
there exists a unique equivariant homomorphism $ [[l]]: \mathcal{T}A 
\rightarrow B $ such that $ [[l]] \sigma_A = l $. 
\end{prop}
An important ingredient in the Cuntz-Quillen approach to cyclic homology \cite{CQ1}, \cite{CQ2}, \cite{CQ3} is the concept of a 
quasifree pro-algebra. The same is true in the equivariant theory. 
\begin{definition} A pro-$ H $-algebra $ R $ is called $ H $-equivariantly quasifree 
if there exists an equivariant splitting homomorphism $ R \rightarrow \mathcal{T}R $ 
for the projection $ \tau_R $.
\end{definition}
We state some equivalent descriptions of equivariantly quasifree pro-$ H $-algebras. 
\begin{theorem} \label{qf}
Let $ H $ be a bornological quantum group and let $ R $ be a pro-$ H $-algebra. The
following conditions are equivalent:
\begin{bnum}
\item[a)] $ R $ is $ H $-equivariantly quasifree.
\item[b)] There exists an equivariant pro-linear map
$ \nabla: \Omega^1(R) \rightarrow \Omega^2(R) $ satisfying
\begin{equation*}
\nabla(a \omega) = a \nabla(\omega), \qquad
\nabla(\omega a) = \nabla(\omega) a - \omega da
\end{equation*}
for all $ a \in R $ and $ \omega \in \Omega^1(R) $.
\item[c)] There exists a projective resolution $ 0 \rightarrow P_1 \rightarrow P_0 \rightarrow R^+ $ 
of the $ R $-bimodule $ R^+ $ of length $ 1 $ in $ \pro(H \LSMod) $.
\end{bnum}
\end{theorem}
A map $ \nabla: \Omega^1(R) \rightarrow \Omega^2(R) $ satisfying condition b) in theorem \ref{qf} is also 
called an equivariant graded connection on $ \Omega^1(R) $. \\
We have the following basic examples of quasifree pro-$ H $-algebras. 
\begin{prop}\label{TAQuasifree} Let $ A $ be any $ H $-algebra. The periodic tensor algebra 
$ \mathcal{T}A $ is $ H $-equivariantly quasifree. 
\end{prop}
An important result in theory of Cuntz and Quillen relates the $ X $-complex of the
periodic tensor algebra $ \mathcal{T}A $ to the standard complex of $ A $ constructed using noncommutative 
differential forms. The comparison between the equivariant $ X $-complex and equivariant differential forms is carried out in 
the same way as in the group case \cite{Voigtepch}. 
\begin{prop}\label{Xdiff}
There is a natural isomorphism $ X_H(\mathcal{T}A) \cong \theta\Omega_H(A) $ such that 
the differentials of the equivariant $ X $-complex correspond to
\begin{alignat*}{2}
\partial_1 &= b - (\id + \kappa)d \qquad &&\text{on}\;\; \theta\Omega_H^{odd}(A)\\
\partial_0 &= - \sum_{j = 0}^{n - 1}
\kappa^{2j} b + B \qquad && \text{on}\;\; \Omega_H^{2n}(A).
\end{alignat*}
\end{prop}
\begin{theorem} \label{homotopyeq}
Let $ H $ be a bornological quantum group and let $ A $ be an $ H $-algebra. 
Then the paracomplexes $ \theta \Omega_H(A) $ and 
$ X_H(\mathcal{T}A) $ are homotopy equivalent. 
\end{theorem}
For the proof of theorem \ref{homotopyeq} it suffices to observe that the corresponding arguments 
in \cite{Voigtepch} are based on the relations obtained in proposition 
\ref{diffformformel}. 

\section{Equivariant periodic cyclic homology} \label{secech}

In this section we define equivariant periodic cyclic homology for bornological quantum groups. 
\begin{definition} Let $ H $ be a bornological quantum group and let 
$ A $ and $ B $ be $ H $-algebras. 
The equivariant periodic cyclic homology of $ A $ and $ B $ is 
\begin{equation*}
HP^H_*(A,B) = 
H_*(\Hom_{\AYD(H)}(X_H(\mathcal{T}(A \rtimes H \rtimes \hat{H})), 
X_H(\mathcal{T}(B \rtimes H \rtimes \hat{H}))).
\end{equation*}
\end{definition}
We write $ \Hom_{\AYD(H)} $ for the space of $ \ayd $-maps and consider the usual 
differential for a $ \Hom $-complex in this definition. Using proposition \ref{covparaadd} it is 
straightforward to check that this yields indeed 
a complex. Remark that both entries in the above $ \Hom $-complex are only paracomplexes. \\
Let us consider the special case that $ H = \D(G) $ is the smooth group algebra of a locally compact group $ G $. 
In this situation the definition of $ HP^H_* $ reduces to the definition of $ HP^G_* $ given in 
\cite{Voigtepch}. This is easily seen using the Takesaki-Takai isomorphism obtained in proposition \ref{TakTak} 
and the results from \cite{Voigtbqg}. \\
As in the group case $ HP^H_* $ is a bifunctor, contravariant in the first variable and
covariant in the second variable. We define $ HP^H_*(A) =
HP^H_*(\mathbb{C}, A) $ to be the equivariant periodic cyclic
homology of $ A $ and $ HP_H^*(A) = HP^H_*(A, \mathbb{C}) $ to be
equivariant periodic cyclic cohomology. 
There is a natural associative product
\begin{equation*}
HP^H_i(A,B) \times HP^H_j(B,C) \rightarrow HP^H_{i + j}(A,C),
\qquad (x,y) \mapsto x\cdot y
\end{equation*}
induced by the composition of maps. Every equivariant homomorphism $ f: A \rightarrow B $
defines an element in $ HP^H_0(A,B) $ denoted by $ [f] $. The
element $ [\id] \in HP^H_0(A,A) $ is denoted $ 1 $ or $ 1_A $. An element $ x \in HP^H_*(A,B) $ is 
called invertible if there
exists an element $ y \in HP^H_*(B,A) $ such that $ x \cdot y =
1_A $ and $ y \cdot x = 1_B $. An invertible element of degree
zero is called an $ HP^H $-equivalence. Such an element
induces isomorphisms $ HP^H_*(A,D) \cong HP^H_*(B,D) $ and $ HP^H_*(D,A)
\cong HP^H_*(D,B) $ for all $ H $-algebras $ D $. 

\section{Homotopy invariance, stability and excision} \label{sechomstabex}

In this section we show that equivariant periodic cyclic homology is homotopy invariant, stable and 
satisfies excision in both variables. Since the arguments carry over from the group case with 
minor modifications most of the proofs will only be sketched. More details can be found in \cite{Voigtepch}. \\
We begin with homotopy invariance. Let $ B $ be a pro-$ H $-algebra and consider the Fr\'echet algebra $ C^\infty[0,1] $ 
of smooth functions on the interval $ [0,1] $. We denote by $ B[0,1] $ the 
pro-$ H $-algebra $ B \hat{\otimes} C^\infty[0,1] $ where the action on 
$ C^\infty[0,1] $ is trivial. A (smooth) equivariant homotopy is an equivariant 
homomorphism $ \Phi: A \rightarrow B[0,1] $ of $ H $-algebras. 
Evaluation at the point $ t \in [0,1] $ yields an equivariant homomorphism 
$ \Phi_t: A \rightarrow B $. Two equivariant homomorphisms from $ A $ to $ B $ are called 
equivariantly homotopic if they can be connected by an equivariant homotopy. 
\begin{theorem}[Homotopy invariance] \label{homotopyinv} Let $ A $ and $ B $ be 
$ H $-algebras and let $ \Phi: A \rightarrow B[0,1] $ be a smooth equivariant homotopy.
Then the elements $ [\Phi_0] $ and $ [\Phi_1] $ in $ HP^H_0(A,B) $ are
equal. Hence the functor $ HP^H_* $ is homotopy invariant in both variables
with respect to smooth equivariant homotopies. 
\end{theorem}
Recall that $ \theta^2 \Omega_H(A) $ is the 
paracomplex $ \Omega^0_H(A) \oplus \Omega^1_H(A) \oplus \Omega^2_H(A)/b(\Omega^3_H(A)) $ with the usual differential $ B + b $ 
and the grading into even and odd forms for any pro-$ H $-algebra $ A $. 
There is a natural chain map $ \xi^2: \theta^2 \Omega_H(A) \rightarrow X_H(A) $. 
\begin{prop} \label{homotopyinv1a} Let $ A $ be an equivariantly quasifree pro-$ H $-algebra. 
Then the map $ \xi^2: \theta^2 \Omega_H(A) \rightarrow X_H(A) $ is a homotopy equivalence. 
\end{prop}
A homotopy inverse is constructed using an equivariant connection for $ \Omega^1(A) $. \\
Now let $ \Phi: A \rightarrow B[0,1] $ be an equivariant homotopy. 
The derivative of $ \Phi $ is an equivariant linear map 
$ \Phi': A \rightarrow B[0,1] $. If we view $ B[0,1] $ as a bimodule over itself the map 
$ \Phi' $ is a derivation with respect to $ \Phi $ in the sense that 
$ \Phi'(ab) = \Phi'(a) \Phi(b) + \Phi(a) \Phi'(b) $ for $ a, b \in A $. 
We define an $ \ayd $-map $ \eta: \Omega^n_H(A) \rightarrow \Omega^{n - 1}_H(B) $ 
for $ n > 0 $ by 
\begin{equation*}
\eta(x \otimes a_0da_1 \dots da_n) = 
\int_0^1 x \otimes \Phi_t(a_0) \Phi'_t(a_1) d\Phi_t(a_2) \cdots d\Phi_t(a_n) dt
\end{equation*}
and an explicit calculation yields the following result. 
\begin{lemma} \label{homotopyinv2} We have $ X_H(\Phi_1) \xi^2 - X_H(\Phi_0) \xi^2 = \partial \eta + \eta \partial $. Hence the chain maps 
$ X_H(\Phi_t) \xi^2: \theta^2 \Omega_H(A) \rightarrow 
X_H(B) $ for $ t = 0,1 $ are homotopic. 
\end{lemma}
Using the map $ \Phi $ we obtain an equivariant homotopy 
$ A \cotimes \mathcal{K}_H \rightarrow (B \cotimes \mathcal{K}_H)[0,1] $
which induces an equivariant homomorphism 
$ \mathcal{T}(A \cotimes \mathcal{K}_H) 
\rightarrow \mathcal{T}((B \cotimes \mathcal{K}_H)[0,1]) $. 
Together with the obvious homomorphism $ \mathcal{T}((B \cotimes \mathcal{K}_H)[0,1]) 
\rightarrow \mathcal{T}(B \cotimes \mathcal{K}_H)[0,1] $ this yields an equivariant homotopy 
$ \Psi: \mathcal{T}(A \cotimes \mathcal{K}_H) \rightarrow \mathcal{T}(B \cotimes \mathcal{K}_H)[0,1] $. 
Since $ \mathcal{T}(A \cotimes \mathcal{K}_H) $ is equivariantly quasifree 
we can apply proposition \ref{homotopyinv1a} and lemma \ref{homotopyinv2} 
to obtain $ [\Phi_0] = [\Phi_1] \in HP^H_0(A,B) $. This finishes the proof of theorem 
\ref{homotopyinv}. \\
Homotopy invariance has several important consequences. Let us call an extension 
$ 0 \rightarrow J \rightarrow R \rightarrow A \rightarrow 0 $ of pro-$ H $-algebras 
with equivariant pro-linear splitting a universal locally nilpotent extension of $ A $ if 
$ J $ is locally nilpotent and $ R $ is equivariantly quasifree. In particular, 
$ 0 \rightarrow \mathcal{J}A \rightarrow \mathcal{T}A \rightarrow A \rightarrow 0 $ is a 
universal locally nilpotent extension of $ A $. Using homotopy invariance one shows that $ HP^H_* $ can be computed 
using arbitrary universal locally nilpotent extensions. \\
Let us next study stability. One has to be slightly careful to formulate correctly the statement of 
the stability theorem since the tensor product of two $ H $-algebras is no longer an $ H $-algebra 
in general. \\  
Let $ H $ be a bornological quantum group and assume that we are given an essential $ H $-module $ V $ 
together with an equivariant bilinear pairing 
$ b: V \times V \rightarrow \mathbb{C} $. Moreover let $ A $ be an $ H $-algebra. Recall from section \ref{secaction} that 
$ l(b;A) = V \cotimes A \cotimes V $ is the $ H $-algebra with multiplication 
$$
(v_1 \otimes a_1 \otimes w_1)(v_2 \otimes a_2 \otimes w_2) 
= b(w_1, v_2)\, v_1 \otimes a_1a_2 \otimes w_2 
$$ 
and the diagonal $ H $-action. We call the pairing $ b $
admissible if there exists an $ H $-invariant vector $ u \in V $ such that $ b(u,u) = 1 $. 
In this case the map $ \iota_A: A \rightarrow l(b; A) $ given by $ \iota(a)
= u \otimes a \otimes u $ is an equivariant homomorphism. 
\begin{theorem}\label{SSStabLemma} Let $ H $ be a bornological quantum group and let $ A $ be an 
$ H $-algebra. For every admissible equivariant bilinear pairing 
$ b: V \times V \rightarrow \mathbb{C} $ the map $ \iota: A \rightarrow 
l(b;A) $ induces an invertible element $ [\iota_A] \in
H_0(\Hom_{\AYD(H)}(X_H(\mathcal{T}A), X_H(\mathcal{T}(l(b;A)))) $.
\end{theorem}
\proof The canonical map $ l(b; A) \rightarrow l(b; \mathcal{T}A) $
induces an equivariant homomorphism $ \lambda_A: \mathcal{T}(l(b; A)) 
\rightarrow l(b; \mathcal{T}A) $. Define the
map $ tr_A: X_H(l(b; \mathcal{T}A)) \rightarrow
X_H(\mathcal{T}A) $ by
\begin{equation*}
tr_A(x \otimes (v_0 \otimes a_0 \otimes w_0)) = b(S^{-1}(x_{(1)}) \cdot w_0,v_0)\, 
x_{(2)} \otimes a_0
\end{equation*}
and
\begin{equation*}
tr_A(x \otimes (v_0 \otimes a_0 \otimes w_0) \, d(v_1 \otimes a_1 \otimes w_1)) = 
b((S^{-1}(x_{(1)}) \cdot w_1, v_0) b(w_0,v_1)\, x_{(2)} \otimes a_0 da_1.
\end{equation*}
In these formulas we implicitly use the twisted trace $ tr_x: l(b) \rightarrow \mathbb{C} $ 
for $ x \in H $ defined by $ tr_x(v \otimes w) = b((S^{-1}(x) \cdot w,v) $. 
The twisted trace satisfies 
the relation 
\begin{equation*}
tr_x(T_0 T_1) = tr_{x_{(2)}}((S^{-1}(x_{(1)}) \cdot T_1) T_0)
\end{equation*}
for all $ T_0,T_1 \in l(b) $. 
Using this relation one verifies that $ tr_A $ defines a chain map. It is clear that $ tr_A $ is $ \hat{H} $-linear and 
it is straightforward to check that $ tr_A $ is $ H $-linear. Let us define $ t_A = tr_A \circ
X_H(\lambda_A) $ and show that $ [t_A] $ is an inverse for 
$ [\iota_A] $. Using the fact that $ u $ is $ H $-invariant one computes 
$ [\iota_A] \cdot [t_A] = 1 $. We have to prove $ [t_A] \cdot
[\iota_A] = 1 $.
Consider the following equivariant homomorphisms 
$ l(b; A) \rightarrow l(b; l(b; A)) $ given by
\begin{align*}
&i_1(v \otimes a \otimes w) = u \otimes v \otimes a \otimes w \otimes u \\
&i_2(v \otimes a \otimes w) = v \otimes u \otimes a \otimes u \otimes w
\end{align*}
As above we see $ [i_1] \cdot [t_{l(b; A)}] = 1 $ and we determine $ [i_2] \cdot [t_{l(b; A)}] =
[t_A] \cdot [\iota_A] $. Let $ h_t $ be the linear map from $ l(b;A) $ into $ l(b; l(b;A)) $ given by 
\begin{align*}
&h_t(v \otimes a \otimes w) = \cos(\pi t/2)^2 u \otimes v \otimes a \otimes w \otimes u 
+ \sin(\pi t/2)^2 v \otimes u \otimes a \otimes u \otimes w \\
&\,- i\cos(\pi t/2) \sin(\pi t/2) u \otimes v \otimes a \otimes u \otimes w 
+ i\sin(\pi t/2) \cos(\pi t/2) v \otimes u \otimes a \otimes w \otimes u 
\end{align*}
The family $ h_t $ depends smoothly on $ t $ and
we have $ h_0 = i_1 $ and $ h_1 = i_2 $. 
Since $ u $ is invariant the map $ h_t $ is in fact equivariant 
and one checks that $ h_t $ is a homomorphism. 
Hence we have indeed defined a smooth homotopy
between $ i_1 $ and $ i_2 $. This
yields $ [i_1] = [i_2] $ and hence $ [t_A] \cdot [\iota_A] = 1 $. \qed \\
We derive the following general stability theorem.
\begin{prop}[Stability] \label{stability} Let $ H $ be a bornological quantum group and let $ A $ be an $ H $-algebra. Moreover let 
$ V $ be an essential $ H $-module and let $ b: V \times V \rightarrow \mathbb{C} $ be a nonzero equivariant bilinear pairing. 
Then there exists an invertible element in $ HP^G_0(A, l(b; A)) $. Hence there are natural isomorphisms
\begin{equation*}
HP^H_*(l(b; A), B) \cong HP^H_*(A,B) \qquad HP^H_*(A,B)
\cong HP^H_*(A, l(b; B))
\end{equation*}
for all $ H $-algebras $ A $ and $ B $.
\end{prop}
\proof Let us write $ \beta: H \times H \rightarrow \mathbb{C} $ for the canonical equivariant bilinear pairing 
introduced in section \ref{secaction}. Moreover we denote by $ b_\tau $ the pairing 
$ b: V_\tau \times V_\tau \rightarrow \mathbb{C} $ where $ V_\tau $ is the space $ V $ 
equipped with the trivial $ H $-action. 
We have an equivariant isomorphism $ \gamma: l(b_\tau; l(\beta; A)) \cong l(\beta; l(b; A)) $ 
given by 
$$
\gamma(v \otimes (x \otimes a \otimes y) \otimes w) = x_{(1)} \otimes x_{(2)} \cdot v \otimes a \otimes y_{(1)} \cdot w \otimes y_{(2)}
$$ 
and using $ \beta(x,y) = \psi(S(y)x) $ as well as the fact that $ \psi $ is right invariant one checks that
$ \gamma $ is an algebra homomorphism. Now we can apply theorem \ref{SSStabLemma} to obtain the assertion. \qed \\
We deduce a simpler description of equivariant periodic cyclic homology in certain cases. A bornological quantum group $ H $ is said to be of 
compact type if the dual algebra $ \hat{H} $ is unital. Moreover let us call $ H $ of semisimple type if 
it is of compact type and the value of the integral for $ \hat{H} $ on $ 1 \in \hat{H} $ is nonzero. For instance, the dual of a 
cosemisimple Hopf algebra $ \hat{H} $ is of semisimple type. 
\begin{prop} \label{defcomp} 
Let $ H $ be a bornological quantum group of semisimple type. Then we have
\begin{equation*}
HP^H_*(A,B) \cong H_*(\Hom_{\AYD(H)}(X_H(\mathcal{T}A), X_H(\mathcal{T}B)))
\end{equation*}
for all $ H $-algebras $ A $ and $ B $.
\end{prop}
\proof Under the above assumptions the canonical bilinear 
pairing $ \beta: \hat{H} \times \hat{H} \rightarrow \mathbb{C} $ is admissible since the element $ 1 \in \hat{H} $ is invariant. \qed \\
Finally we discuss excision in equivariant periodic cyclic homology. Consider an extension 
\begin{equation*}
   \xymatrix{
     K\;\; \ar@{>->}[r]^{\iota} & E \ar@{->>}[r]^{\pi} & Q 
     }
\end{equation*}
of $ H $-algebras equipped with an equivariant linear splitting $ \sigma: Q \rightarrow E $ 
for the quotient map $ \pi: E \rightarrow Q $. \\
Let $ X_H(\mathcal{T}E:\mathcal{T}Q) $ be the kernel of the map $ X_H(\mathcal{T}\pi):
X_H(\mathcal{T}E) \rightarrow X_G(\mathcal{T}Q) $ induced by $ \pi $. The splitting $ \sigma $ 
yields a direct sum decomposition 
$ X_H(\mathcal{T}E) = X_H(\mathcal{T}E:\mathcal{T}Q) \oplus X_H(\mathcal{T}Q) $ 
of $ \ayd $-modules. The resulting extension 
$$ 
   \xymatrix{
      X_H(\mathcal{T}E:\mathcal{T}Q)\;\; \ar@{>->}[r] & X_H(\mathcal{T}E) \ar@{->>}[r] & X_H(\mathcal{T}Q) 
     }
$$ 
of paracomplexes induces long exact sequences in homology in both variables. 
There is a natural covariant map 
$ \rho: X_H(\mathcal{T}K) \rightarrow X_H(\mathcal{T}E:\mathcal{T}Q) $ of paracomplexes and we 
have the following generalized excision theorem. 
\begin{theorem}\label{Excision2} The map $ \rho: X_H(\mathcal{T}K) \rightarrow
X_H(\mathcal{T}E:\mathcal{T}Q) $ is a homotopy equivalence.
\end{theorem}
This result implies excision in equivariant periodic cyclic homology.
\begin{theorem}[Excision]\label{Excision} Let $ A $ be an $ H $-algebra and let 
$ (\iota, \pi): 0 \rightarrow K \rightarrow E \rightarrow Q \rightarrow 0 $ be
an extension of $ H $-algebras with a bounded linear splitting. Then there are two natural exact sequences
\begin{equation*}
\xymatrix{
 {HP^H_0(A,K)\;} \ar@{->}[r] \ar@{<-}[d] &
      HP^H_0(A,E) \ar@{->}[r] &
        HP^H_0(A,Q) \ar@{->}[d] \\
   {HP^H_1(A,Q)\;} \ar@{<-}[r] &
    {HP^H_1(A,E)}  \ar@{<-}[r] &
     {HP^H_1(A,K)} \\
}
\end{equation*}
and
\begin{equation*}
\xymatrix{
    {HP^H_0(Q,A)\;} \ar@{->}[r] \ar@{<-}[d] &
       HP^H_0(E,A) \ar@{->}[r] &
          HP^H_0(K,A) \ar@{->}[d] \\
    {HP^H_1(K,A)\;} \ar@{<-}[r] &
      {HP^H_1(E,A)}  \ar@{<-}[r] &
        {HP^H_1(Q,A)} \\
}
\end{equation*}
The horizontal maps in these diagrams are induced by the maps
in the extension.
\end{theorem} 
In theorem \ref{Excision} we only require a linear splitting for the quotient 
homomorphism $ \pi: E \rightarrow Q $. Taking double crossed products of the extension given in theorem \ref{Excision} yields 
an extension 
\begin{equation*}
   \xymatrix{
     K \cotimes \mathcal{K}_H \;\; \ar@{>->}[r] & E \cotimes \mathcal{K}_H \ar@{->>}[r] & Q \cotimes \mathcal{K}_H 
     }
\end{equation*}
of $ H $-algebras with a linear splitting. Using lemma 
\ref{omegaprojective} one obtains an equivariant linear splitting for this extension. 
Now we can apply theorem \ref{Excision2} to this extension and obtain 
the claim by considering long exact sequences in homology. \\
For the proof of theorem \ref{Excision2} one considers the left ideal $ \mathfrak{L} \subset \mathcal{T}E $ generated by
$ K \subset \mathcal{T}E $. The natural projection 
$ \tau_E: \mathcal{T}E \rightarrow E $ induces an equivariant 
homomorphism $ \tau: \mathfrak{L} \rightarrow K $ and one obtains an extension 
\begin{equation*}
  \xymatrix{
     N\;\; \ar@{>->}[r] & \mathfrak{L} \ar@{->>}[r]^\tau & K 
     }
\end{equation*}
of pro-$ H $-algebras. The pro-$ H $-algebra $ N $ 
is locally nilpotent and theorem \ref{Excision2} follows from the following assertion. 
\begin{theorem}\label{Excision3} With the notations as above we have
\begin{bnum}
\item[a)] The pro-$ H $-algebra $ \mathfrak{L} $ is equivariantly quasifree.
\item[b)] The inclusion map $ \mathfrak{L} \rightarrow \mathcal{T}E $ induces a
homotopy equivalence $ \psi: X_H(\mathfrak{L}) \rightarrow
X_H(\mathcal{T}E:\mathcal{T}Q) $.
\end{bnum}
\end{theorem}
In order to prove the first part of theorem \ref{Excision3} one constructs explicitly a 
projective resolution of $ \mathfrak{L} $ of length one. 

\bibliographystyle{plain}

\end{document}